\documentclass[11pt,twoside,reqno]{amsart}
\usepackage[colorlinks=true,citecolor=blue]{hyperref}
\usepackage{mathptmx, amsmath, amssymb, amsfonts, amsthm, mathptmx, enumerate, color}
\usepackage{graphicx}
\usepackage{epstopdf}

\newtheorem{theorem}{Theorem}[section]

\newtheorem{proposition}[theorem]{Proposition}
\newtheorem{corollary}[theorem]{Corollary}

\theoremstyle{definition}
\newtheorem{definition}[theorem]{Definition}

\newtheorem{example}[theorem]{Example}

\newtheorem{remark}[theorem]{Remark}
\numberwithin{equation}{section}

\begin{document}
\setcounter{page}{1}

\vspace*{2.0cm}
\title[The fixed point property for $(c)$-mappings and unbounded sets]
{The fixed point property for $(c)$-mappings and unbounded sets}
\author[S. Atailia, A. Dehici, N. Redjel]{ Sami Atailia$^{1,*}$, Abdelkader Dehici$^{2}$, Najeh Redjel$^2$}
\maketitle
\vspace*{-0.6cm}

\begin{center}
{\footnotesize

$^1$Department of Mathematics, University of Boumerdes, 35000 Boumerdes, Algeria.
\\
$^2$Laboratory of Informatics and Mathematics,
University of Souk-Ahras, P.O.Box 1553, Souk-Ahras, 41000, Algeria.

}\end{center}


\vskip 4mm {\footnotesize \noindent {\bf Abstract.}
We prove that a closed convex subset $C$ of a real Hilbert space $X$ has the fixed point property for $(c)$-mappings if and only if $C$ is bounded. Some convergence results about the iterations are obtained.


 }

\renewcommand{\thefootnote}{}
\footnotetext{ $^*$Corresponding author.
\par
E-mail addresses: s.atailia@univ-boumerdes.dz (S. Atailia), dehicikader@yahoo.fr  (A. Dehici), n.radjal@univ-soukahras.dz (N. Redjel).
\par
Received January ... ; Accepted ... }

\section{Introduction}
Let $X$ be a real Banach space with norm $\|.\|$ and let $C$ be a nonempty subset of $X$. A mapping $T: C \longrightarrow C$ is said to be nonexpansive if $\| Tx - Ty \| \leq \| x - y \|$ for all $x, y \in C$. $T$ is said to be  a $(c)$-mapping if there exist $a, c \in [0, 1], c > 0$ and $a + 2c = 1$ such that
\begin{equation}\label{eq 1.1}
\|Tx - Ty \| \leq a \| x - y \| + c (\|Tx - y \| + \|Ty - x \|).
\end{equation}
\noindent for all $x, y \in C$.
\vskip 0.3 cm
\noindent A closed convex subset $C$ of $X$ is said to have the fixed point property for nonexpansive mappings (in short, FPP) if every nonexpansive mapping $T: C \longrightarrow C$ has at least a fixed point in $C$ (see \cite{Bro, Goe2, Goe3, Goh, Kha1, Kir1, Rei4}).
\vskip 0.3 cm
\noindent A closed convex subset $C$ of $X$ is said to have the fixed point property for $(c)$-mappings (in short, $(c)$-FPP) if every $(c)$-mapping $T: C \longrightarrow C$ has at least a fixed point in $C$ (see \cite{Rei6, Smy}).
\vskip 0.3 cm
\noindent It is an open problem whether these two fixed point properties hold simultaneously. The answer is affirmative if $C$ is a bounded set of uniformly convex Banach space (see \cite{Goe3, Kha1, Kir1}). However, the situation seems to be unknown when $C$ is unbounded. The contributions related to this subject are very few. In 1980, W. Ray (see \cite{Ray}) proved that the boundedness of $C$ characterizes FPP in Hilbert spaces. Ray's result was simplified by R. Sine \cite{Sin} who observed that the metric projection in Hilbert space is nonexpansive.  After that, T. Benavides (see \cite{Ben}) established the same result in the Banach space $c_0$.
In \cite{Smy}, M. A. Smyth investigated the existence of fixed points for $(c)$-mappings defined on weakly compact convex subsets which does not have necessarily normal structure and he wondered about assumptions on $C$ to be a weakly compact convex subset having FPP, to posses $(c)$-FPP. Recall that in the setting of Banach space $L^{1}([0, 1])$, the weakly compact convex subset
\vskip 0.3 cm
\centerline {$ C = \{ f \in L^{1}([0, 1]): 0 \leq f \leq 2, a.e, \displaystyle \int_0^{1} f(t) dt = 1\}$}
\vskip 0.3 cm
\noindent fails to have FPP (see \cite{Als}) but we do not know if $C$ has $(c)$-FPP.
\vskip 0.3 cm
\noindent In this paper, by the insights in the contributions of W. Takahashi et al in \cite{Tak}, we prove the variant of Ray's result for $(c)$-mappings. Some convergence of iterations associated to $(c)$-mappings are studied. Finally, we conclude this work by asking some interesting questions.
\vskip 0.3 cm
\section{Main Results}
\noindent First of all, let us define the concept of firmly nonexpansive mappings.
\vskip 0.3 cm
\begin{definition} \label{def 1.1}\rm Let $C$ be a nonempty subset of a Banach space $X$. A mapping $T: C \longrightarrow C$ is said to be $\lambda$-firmly nonexpansive $(\lambda \in (0, 1)$ if
	\begin{equation}\label{eq 2.1}
	\|Tx- Ty\|\leq \| (1 - \lambda) (x - y) + \lambda (Tx - Ty)\|
	\end{equation}
	\noindent for all $x, y \in C$. $T$ is said to be firmly nonexpansive if $T$ is $\lambda$-firmly nonexpansive for all $\lambda \in (0, 1)$.
\end{definition}

\vskip 0.3 cm
\begin{remark} \label{rem 1.1}\rm It is obvious that a $\lambda$-firmly nonexpansive mapping is nonexpansive while the converse is in general not true (it suffices to take $T: X \longrightarrow X$ defined by $Tx = -x$).
	\vskip 0.3 cm
	\noindent For more details on $\lambda$-firmly nonexpansive mappings, we quote \cite{Ari, Goe2, Goe3, Rei2, Sma}).
	\vskip 0.3 cm
	\noindent The first result in this section is the following proposition.
\end{remark}

\vskip 0.3 cm
\begin{proposition} \label{prop 1.1}\rm Let $C$ be a nonempty subset of a Banach space $X$. Then every $\lambda$-firmly nonexpansive mapping is a $(c)$-mapping.
\end{proposition}
\proof Let $x, y \in C$. Since $T$ is a $\lambda$-firmly nonexpansive mapping, we have
\begin{align}
\|Tx - Ty\| \leq & \|\lambda (Tx - Ty) + (1 - \lambda) (x - y)\| \nonumber \\
= & \|[(1 - \lambda) x + \lambda Tx] - y + \lambda y - \lambda Ty]\|\nonumber \\
= & \|(1 - \lambda) [(1 - \lambda)x + \lambda Tx- y] + \lambda [(1 - \lambda) x + \lambda Tx - Ty]\|\nonumber \\
\leq & (1 - \lambda) \|[(1 - \lambda)x + \lambda Tx] - y\|  + \lambda \|(1 - \lambda) x + \lambda Tx - Ty\|\nonumber \\
\leq & (1 - \lambda) \|[(1 - \lambda)x + \lambda Tx] - [(1 - \lambda)y + \lambda y]\| \nonumber \\
& + \lambda \|[(1 - \lambda) x + \lambda Tx] - [(1 - \lambda)Ty + \lambda Ty]\|\nonumber \\
\leq & (1 - \lambda) [(1 - \lambda)\| x - y\| + \lambda \|Tx - y\|] \nonumber \\
& + \lambda [(1 - \lambda)\|x - Ty\|  + \lambda \|Tx - Ty\|] \nonumber \\
=& (1 - \lambda)^{2}\| x - y\| + \lambda (1 - \lambda) \|Tx - y\| + \lambda (1 - \lambda)\|x - Ty\| \nonumber \\
& + \lambda^{2} \|Tx - Ty\|. \nonumber
\end{align}
\noindent So,
\vskip 0.3 cm
\centerline{$  (1 - \lambda^{2}) \| Tx - Ty\| \leq (1 - \lambda)^{2} \|x - y\| + \lambda (1 - \lambda)\|Tx - y\| + \lambda (1 - \lambda) \|x - Ty\|.$}
\vskip 0.3 cm
\noindent Therefore
\vskip 0.3 cm
\begin{equation}\label{eq 2.2}
\| Tx - Ty\| \leq \displaystyle \frac{1 - \lambda}{ 1 + \lambda} \|x - y\| + \displaystyle  \frac{\lambda}{ 1 + \lambda} (\|Tx - y\| +  \|x - Ty\|).
\end{equation}
\vskip 0.3 cm
\noindent which means that $T$ is a $(c)$-mapping.
\vskip 0.3 cm
\noindent The following example shows that the class of $(c)$-mappings is wider than that of firmly nonexpansive mappings.
\vskip 0.3 cm

\noindent \vskip 0.3 cm
\begin{example} \label{exa 1.1}\rm (see \cite{Suz}) Let $(X, \|.\|)= (\mathbb{R}, |.|)$ and $C= [0, 3]$. Define $T: [0, 3] \longrightarrow [0, 3]$ by
	\vskip 0.3 cm
	\centerline{$Tx= \left\{
		\begin{array}{ll}
		0 & \hbox{if} \ x \in [0,3[, \\
		1 & \hbox{if} \ x= 3.
		\end{array}
		\right.$}
	\vskip 0.3 cm
	\noindent A simple calculation shows that $T$ is a $(c)$-mapping for $c = \displaystyle \frac{1}{2}$ and $c = \displaystyle \frac{1}{3}$. However, $T$ is not firmly nonexpansive, since $T$ is not nonexpansive.
\end{example}
\vskip 0.3 cm
\noindent The following theorem was established by W. Takahashi et al in \cite{Tak}.

\vskip 0.3 cm
\begin{theorem} \label{the 1.1}\rm Let $H$ be a Hilbert space and let $C$ be a nonempty closed convex subset of $H$. Then the following conditions are equivalent.
	\vskip 0.3 cm
	\noindent $(\imath)$ Every firmly nonexpansive mapping $T: C \longrightarrow C$ has a fixed point in $C$.
	\vskip 0.3 cm
	\noindent $(\imath \imath)$ $C$ is bounded.
\end{theorem}

\vskip 0.3 cm
\noindent The main theorem of this paper is the following.

\vskip 0.3 cm
\begin{theorem} \label{the 1.1}\rm Let $H$ be a Hilbert space and let $C$ be a nonempty closed convex subset of $H$. Then the following conditions are equivalent.
	\vskip 0.3 cm
	\noindent $(\imath)$ Every $(c)$-mapping $T: C \longrightarrow C$ has a fixed point in $C$.
	\vskip 0.3 cm
	\noindent $(\imath \imath)$ $C$ is bounded.
\end{theorem}
\proof $(\imath \imath) \Longrightarrow (\imath)$ Since $C$ is bounded then $C$ is a weakly compact convex subset of $H$ (which has a normal structure). So the result is an immediate consequence of Theorem 2 in \cite{Bog} (see also \cite{Kir1}).
\vskip 0.3 cm
\noindent $(\imath) \Longrightarrow (\imath \imath)$ Assume that $C$ is unbounded. By using Theorem 2.5, there exists a free fixed point firmly nonexpansive mapping $T: C \longrightarrow C$. But following Proposition 2.3, $T$ is a $(c)$-mapping which contradicts $(\imath)$. Hence $C$ must be bounded.

\vskip 0.3 cm
\begin{corollary} \label{cor 1.1}\rm Let $C$ be a nonempty closed convex subset of a Hilbert space $H$ and let $T: C \longrightarrow C$ be a mapping satisfying
	\vskip 0.3 cm
	\begin{equation}\label{eq 2.3}
	\ 9 \|Tx - Ty\|^{2}\leq \|x - y \|^{2} + \|Tx - y \|^{2} + \|Ty - x \|^{2}
	\end{equation}
	\vskip 0.3 cm
	\noindent If $C$ is bounded then $T$ has a (unique) fixed point in $C$.
\end{corollary}
\proof For the uniqueness, assume that $T$ has two distincts fixed points $z_1, z_2 \in C$ such that $z_1 \neq z_2$. Then
\vskip 0.3 cm
\centerline {$ 9 \|z_1 - z_2\|^{2}\leq \|z_1 - z_2\|^{2} + \|z_1 - z_2 \|^{2} + \|z_1 - z_2 \|^{2} = 3 \|z_1 - z_2 \|^{2}.$}
\vskip 0.3 cm
\noindent which is a contradiction.
\vskip 0.3 cm
\noindent Now, if $T$ satisfies $(4)$, then

\vskip 0.3 cm
\centerline {$ \|Tx - Ty\|^{2}\leq \displaystyle \frac{1}{9} (\|x - y \|^{2} + \|Tx - y \|^{2} + \|Ty - x \|^{2})$.}
\vskip 0.3 cm
\noindent Consequently,
\vskip 0.3 cm
\centerline {$ \|Tx - Ty\|\leq \displaystyle \frac{1}{3} (\|x - y \|^{2} + \|Tx - y \|^{2} + \|Ty - x \|^{2})^{\frac{1}{2}}$.}
\vskip 0.3 cm
\noindent By using the inequality
\vskip 0.3 cm
\centerline {$ \sqrt{x^{2} + y^{2} + z^{2}}\leq x + y + z \ \hbox{for all} \ x, y, z \geq 0$,}
\vskip 0.3 cm
\noindent we get
\vskip 0.3 cm
\begin{equation}\label{eq 2.4}
\|Tx - Ty\|\leq \displaystyle \frac{1}{3} (\|x - y \| + \|Tx - y \| + \|Ty - x \|),
\end{equation}
\vskip 0.3 cm
\noindent which proves that $T$ is a $(c)$-mapping. Now, the result is an immediate consequence of the implication $ (\imath \imath) \Longrightarrow (\imath)$ of Theorem 2.6.

\vskip 0.3 cm
\begin{definition} \label{def 1.1}\rm Let $X$ be a Banach space. The modulus $\delta$ of convexity of $X$ is defined by
	\vskip 0.3 cm
	\centerline{$ \delta(\epsilon) = \inf \{ 1 - \displaystyle \frac{\| x + y\|}{2}: \|x\| \leq 1, \|y\| \leq 1, \|x - y \| \geq \epsilon\}$}
	\vskip 0.3 cm
	\noindent for every $0 \leq \epsilon \leq 2$. A Banach space $X$ is said to be uniformly convex if $\delta(\epsilon) > 0$ for all $\epsilon \in (0, 2]$.
	\vskip 0.3 cm
	\noindent For a mapping $T: C \longrightarrow C$, we define the orbit $O(x_0)$ of $x_0 \in C$ by $O(x_0) = \{T^{n}x_0\}_{n \geq 0} \ (T^{0} x_0 = x_0)$.
\end{definition}

\vskip 0.3 cm
\begin{theorem} \label{the 1.1}\rm Let $X$ be a uniformly convex Banach space and let $C$ be a closed convex subset of $X$. Assume that $T: C \longrightarrow C$ is a $(c)$-mapping satisfying
	\vskip 0.3 cm
	\noindent $(\imath)$ There exists $x_0 \in C$ such that $O(x_0)$ is bounded.
	\vskip 0.3 cm
	\noindent Then $T$ has a fixed point in $C$.
	\vskip 0.2 cm
	\noindent If moreover, $C = -C$ and $T$ is an odd mapping satisfying
	\vskip 0.3 cm
	\noindent $(\imath \imath)$ For all integer $i \geq 1$ and all $x, y \in C,$ the sequence $(\|T^{n + i}x - T^{n}y\|)_n$ is decreasing.
	\vskip 0.3 cm
	\noindent Then the Picard sequence $(T^{n}(x_0))_n$ converges in norm to a fixed point of $T$.
\end{theorem}
\proof $(\imath)$ Since $X$ is uniformly convex, then the asymptotic center $A(C, (T^{n}(x_0))_n)$, associated to the Picard sequence $(T^{n}(x_0))_n$, is a singleton (see assertions $(a)$ and $(c)$ of Theorem 5.2 in \cite{Goe3}). On the other hand, Lemma 2 in \cite{Bae} shows that
\vskip 0.3 cm
\begin{equation}\label{eq 2.5}
\displaystyle \lim_{n \longrightarrow + \infty} \|Ty_n - y_n\| = 0.
\end{equation}
\vskip 0.3 cm
\noindent where $y_n = T^{n}(x_0)$. A simple argument implies that there exists $z_0 \in C$ such that $Tz_0 = z_0$ which proves the first claim.
\vskip 0.3 cm
\noindent $(\imath \imath)$ Now, since $C$ is a nonempty convex subset with $C = -C$, we have
\vskip 0.3 cm
\centerline{$0 = \displaystyle \frac{x + (-x)}{2} \in C.$}
\vskip 0.3 cm
\noindent Next, since $T$ is a $(c)$-mapping then for all $x \in C$, we have
\vskip 0.3 cm
\centerline{$\|Tx - T0\| \leq a \| x - 0\| + c(\|Tx - 0\| + \|T0 - x\|)$.}
\vskip 0.3 cm
\noindent But $T$ is odd, so $T0 = 0$. Therefore
\vskip 0.3 cm
\centerline{$\|Tx\| \leq a \| x\| + c\|Tx\| + c \|x\|$.}
\vskip 0.3 cm
\noindent This leads to
\vskip 0.3 cm
\centerline{$\|Tx\| \leq (\displaystyle \frac{a + c}{1 - c}) \|x\|$.}
\vskip 0.3 cm
\begin{equation}\label{eq 2.6}
= \|x\|.
\end{equation}
\vskip 0.3 cm
\noindent Thus, by induction, we deduce that the sequence $(\|T^{n}x\|)_n$ is decreasing in $[0, + \infty[$, and
\vskip 0.3 cm
\begin{equation}\label{eq 2.7}
\displaystyle \lim_{n \longrightarrow + \infty} \|T^{n}x\| = \gamma \geq 0.
\end{equation}
\vskip 0.3 cm
\noindent Furthermore, since $T$ is odd, we have $T^{n}(-x) = - T^{n}x$ for all integer $n \geq 1$. So, by replacing $y$ by $-x$ in $(\imath \imath)$, we observe that the sequence $(\|T^{n + i}x - T^{n}(-x)\| = \|T^{n + i}x + T^{n}x\|)_n$ is non-increasing for all fixed integer $i \geq 1$.
\vskip 0.3 cm
\noindent Afterwards, by the triangle inequality, we infer that
\vskip 0.3 cm
\begin{equation}\label{eq 2.8}
\|T^{n + i}x_0 - T^{n}x_0\|\leq \displaystyle \sum_{k = 1}^{i} \|T^{n + k}x_0 - T^{n + k - 1}x_0\|.
\end{equation}
\vskip 0.3 cm
\noindent It follows, from Lemma 2 in \cite{Bae} that for all fixed integer $i$, we have
\vskip 0.3 cm
\begin{equation}\label{eq 2.9}
\displaystyle \lim_{n \longrightarrow + \infty} \|T^{n + i}x_0 - T^{n}x_0\| = 0.
\end{equation}
\vskip 0.3 cm
\noindent The rest of the proof is similar to that given in Theorem 1.1 in \cite{Bai}.
\vskip 0.3 cm
\noindent The next example shows that the hypothesis of uniform convexity is important in Theorem 2.9.

\vskip 0.3 cm
\begin{example} \label{exa 1.1}\rm Let $X = C([0, 1])$ equipped with the sup norm and let
	\vskip 0.3 cm
	\centerline{$C = \{f \in X : f(0) = 0\}$}
	\vskip 0.3 cm
	\noindent and let $T: C \longrightarrow C$ defined by $Tf(t) = tf(t)$.
	\vskip 0.3 cm
	\noindent Clearly, $C$ is a closed convex subset of $X$ with $C = -C$. In addition, $T$ is an odd $(c)$-mapping (see Example in \cite{Bae}). The formula
	\vskip 0.3 cm
	\begin{equation}\label{eq 2.10}
	T^{n}f(t) = t^{n}f(t) \ (n \geq 1)
	\end{equation}
	\vskip 0.3 cm
	\noindent shows that the orbit $O(f)$ of any $f \in C$ is bounded. On the other hand, it is obvious that $T$ is also nonexpansive then for any fixed integer $i \geq 1$, the sequence $(\|T^{n + i}f_1 - T^{n}f_2\|)_n$ is decreasing and 0 is the unique fixed point of $T$ in $C$. But if we take $f_0(t) = \sin (t\frac{\pi}{2})$ then $f_0 \in C$ and $T^{n}f_0(t) = t^{n} \sin (t\frac{\pi}{2})$ does not converge to 0 in $X$ since $\|T^{n}f_0- 0\| = 1 \nrightarrow 0$.
\end{example}

\vskip 0.3 cm
\begin{definition} \label{def 1.1}\rm A mapping $T: C \longrightarrow C$ is called Chatterjea mapping if $T$ is a $(c)$-mapping with $c = \displaystyle \frac{1}{2}$.
\end{definition}
\vskip 0.3 cm
\noindent In the sequel, we will denote by $R(I - T)$ the range of the mapping $I - T: C \longrightarrow X$.

\vskip 0.3 cm
\noindent Now, we are in a position to state our next result

\vskip 0.3 cm
\begin{theorem} \label{the 1.1}\rm Let $C$ be a closed convex subset of a Banach space $X$ and let $T: C \longrightarrow C$ be a $(c)$-mapping. Then
	\vskip 0.3 cm
	\noindent I) If $X$ is uniformly convex then
	\vskip 0.3 cm
	\centerline{$0 \in R(I - T) \Longleftrightarrow O(x_0)$ is bounded for some $x_0$.}
	\vskip 0.3 cm
	\noindent II) If $T$ is a Chatterjea mapping satisfying the following assumptions:
	\vskip 0.3 cm
	\noindent $\mathcal{H}_1)$ For all $x \in C$ and all integer $k \geq 2$, the sequence $(\|T^{n + k}x - T^{n}x\|)_n$ is decreasing.
	\vskip 0.3 cm
	\noindent $\mathcal{H}_2)$ There exists an integer $k_0 \geq 1$ such that $T^{n}$ is uniformly lipschitzian for all $n \geq k_0$.
	\vskip 0.3 cm
	\noindent $a)$ Then we have
	\vskip 0.3 cm
	\centerline {$0 \notin \overline{R(I - T)} \Longleftrightarrow \displaystyle \lim_{n \longrightarrow + \infty}\frac{\|T^{n}x\|}{n} = \alpha > 0$ for all $ x \in C$.}
	\vskip 0.3 cm
	\noindent $b)$ If $X$ is uniformly convex. Then
	\vskip 0.3 cm
	\centerline{$0 \in \overline{R(I - T)} \  \hbox{and} \ 0 \notin R(I - T) \Longleftrightarrow \displaystyle \lim_{n \longrightarrow + \infty}{\|T^{n}x\|}= \infty \ \hbox{and} \
		\displaystyle \lim_{n \longrightarrow + \infty}\frac{\|T^{n}x\|}{n} = 0$}
	\vskip 0.3 cm
	\noindent  for all $ x \in C$.
	\proof The proof of I) can be obtained by combining Lemma 2 in \cite{Bae} and the equivalence between assertions $(a)$ and $(c)$ of Theorem 5.2 in \cite{Goe3}.
	\vskip 0.3 cm
	\vskip 0.3 cm
	\noindent Now, we will prove II)
	\vskip 0.3 cm
	\noindent II) $ a)\Longrightarrow$  Assume that
	\vskip 0.3 cm
	\centerline {$\displaystyle \lim_{n \longrightarrow + \infty} \frac{\|T^{n}x_0\|}{n} = 0$ for some $x_0 \in C$.}
	\vskip 0.3 cm
	\noindent From Corollary 2.8 in \cite{Deh}, we have
	\vskip 0.3 cm
	\centerline {$\displaystyle \lim_{n \longrightarrow + \infty} \|T^{n + 1}x_0 - T^{n} x_0\| = 0$.}
	\vskip 0.3 cm
	\noindent But
	\vskip 0.3 cm
	\centerline {$\displaystyle \lim_{n \longrightarrow + \infty} \|T(T^{n}x_0) - T^{n} x_0\| = \displaystyle \lim_{n \longrightarrow + \infty} \|(I - T)(T^{n}x_0)\| = 0$}
	\vskip 0.3 cm
	\noindent and
	\vskip 0.3 cm
	\centerline {$(I - T)(T^{n}x_0) \in R(I - T),$}
	\vskip 0.3 cm
	\noindent so
	\vskip 0.3 cm
	\centerline {$ 0 \in \overline{R(I - T)}.$}
	\vskip 0.3 cm
	\noindent II) $ a)\Longleftarrow$  If $0 \in \overline{R(I - T)}.$ Then
	\vskip 0.3 cm
	\centerline {$ 0 = \inf \{ \|y\|: y \in \overline{R(I - T)}\},$}
	\vskip 0.3 cm
	\noindent this implies the existence of a sequence $(x_k)_k$ in $C$ such that
	\vskip 0.3 cm
	\centerline {$\displaystyle \lim_{k \longrightarrow + \infty} \|Tx_k - x_k\|= 0.$}
	\vskip 0.3 cm
	\noindent On the other hand, since
	\vskip 0.3 cm
	\begin{equation}\label{eq 2.11}
	T^{n}x_k = x_k + \displaystyle \sum_{ s = 1}^{n} (T - I) T^{s - 1}x_k,
	\end{equation}
	\vskip 0.3 cm
	\noindent and $T^{n}$ is uniformly lipschitzian for $n \geq k_0$, then
	\vskip 0.3 cm
	\begin{equation}\label{eq 2.12}
	\|T^{n}x\| \leq \|T^{n}x_0\| + M\|x - x_0\| (M > 1),
	\end{equation}
	\vskip 0.3 cm
	\noindent But from (\ref{eq 2.11}), we have
	\vskip 0.3 cm
	\begin{align}
	\|T^{n}x_k\| \leq & \|x_k\| + \displaystyle \sum_{s = 1}^{n} \|(T-I)T^{s-1}x_k\|, \\
	\leq & \|x_k\| + n \|Tx_k - x_k\|,
	\end{align}
	\vskip 0.3 cm
	\noindent it follows that
	\vskip 0.3 cm
	\begin{align}
	\|T^{n}x\| \leq & \|T^{n}x_k\| + M\|x_k - x\|, \\
	\leq & \|x_k\| + n \|Tx_k - x_k\| + M \|x - x_k\|,
	\end{align}
	\vskip 0.3 cm
	\noindent By dividing by $n$, we get
	\begin{align}
	\displaystyle \frac{\|T^{n}x\|}{n} \leq & \displaystyle \frac{\|x_k\|}{n}+ \displaystyle \frac{M\|x_k - x\|}{n} + \|Tx_k-x_k\|.
	\end{align}
	\noindent For a fixed integer $k \geq 1$, letting $n \longrightarrow + \infty$, we infer that
	\vskip 0.3 cm
	\begin{equation}\label{eq 2.17}
	\displaystyle \limsup_n \frac{\|T^{n}x\|}{n}\leq   \|Tx_k - x_k\|,
	\end{equation}
	\vskip 0.3 cm
	\noindent Now, by letting $k \longrightarrow + \infty$, we obtain that
	\vskip 0.3 cm
	\begin{equation}\label{eq 2.18}
	\displaystyle \lim_n \frac{\|T^{n}x\|}{n} = 0.
	\end{equation}
	\vskip 0.3 cm
	\noindent II) $ b)$ Can be deduced immediately from I) and $a)$ of II).
	\vskip 0.3 cm
	\noindent To illustrate Theorem 2.11, we give the following examples
	\vskip 0.3 cm
	\begin{example} \label{exa 1.1}\rm Let $(X, \|.\|) = (\mathbb{R}, |.|)$ and let $T: \mathbb{R}\longrightarrow \mathbb{R}$ defined by $Tx = x + a \ (a \neq 0)$. It is easy to see that $T$ is a free fixed point $(c)$-mapping. Obviously, we have
		\vskip 0.3 cm
		\begin{equation}\label{eq 2.19}
		\displaystyle \lim_n |T^{n}x| = \displaystyle \lim_n |x + na| =  \infty.
		\end{equation}
	\end{example}
	\vskip 0.3 cm
	\noindent which illustrates the assertion I) of Theorem 2.12. This example illustrates also assertion a) of II) since in this case, we have $\overline{R(I - T)} = R(I - T) = \{-a\}$.
	\vskip 0.3 cm
	\begin{example} \label{exa 1.1}\rm Let $X = C([0, 1])$ equipped with the sup norm and
		\vskip 0.3 cm
		\centerline{$ C = \{f \in X: f(0) = 0\leq f(t) \leq f(1) = 1 \}$}
	\end{example}
	\vskip 0.3 cm
	\noindent and let $T: C \longrightarrow C$ be defined as in Example 2.10. Then $T$ is a free fixed point $(c)$-mapping. But since
	\vskip 0.3 cm
	\centerline {$\displaystyle \lim_n \|T^{n + 1}x - T^{n}x\| = 0$ (by Lemma 2 in \cite{Bae}),}
	\vskip 0.3 cm
	\noindent we infer that $0 \in \overline{R(T - I)}$. Next, the fact that all orbits are bounded in this case, we get
	\vskip 0.3 cm
	\centerline{$ \displaystyle \lim_n \displaystyle \frac{\|T^{n}x\|}{n} = 0.$}
	\vskip 0.3 cm
	\noindent This contradicts b) of II). Indeed, $C([0, 1])$ is not uniformly convex.
	\section{Some questions}
	\vskip 0.3 cm
	\noindent We conclude this work by the following interesting questions
	\vskip 0.3 cm
	\noindent {\it Question 1:} Does Banach space $L^{1}([0, 1])$ have $(c)$-FPP?
	\vskip 0.3 cm
	\noindent {\it Question 2 :} Does Benavides's result in $c_0$ hold for $(c)$-mappings?
	\vskip 0.3 cm
	\noindent {\it Question 3 :} Can we extend Benavides's result to orthogonally convex spaces? (see \cite{Smy} for the definition)
	\vskip 0.3 cm
	\noindent {\it Question 4 :} Let $X$ be a Banach space and let $T : X\longrightarrow X$ be a $(c)$-mapping. Is it true that $\overline{R(I - T)}$ is a convex subset of $X$.

	\vskip 0.3 cm
	\noindent Recall that when $T$ is nonexpansive, the convexity of $\overline{R(I - T)}$ was proved by A. Pazy in the case of Hilbert spaces (see \cite{Paz}) and the result was generalized by S. Reich to the setting of uniformly convex Banach spaces (see \cite{Rei3}).

	\vskip 0.3 cm
	\noindent {\it Question 5 :} Let $C$ be a nonempty closed convex subset of a Hilbert space $X$ and let $\mathcal{S}$ be a representation of a semigroup $S$ of $(c)$-mappings on $C$. Suppose that $\{T_sc : s \in S\}$ is relatively weakly compact for some $c \in C$. Does $F(\mathcal{S}) \neq \emptyset$? (Here $F(\mathcal{S})$ is the set of common fixed points of $\mathcal{S}$).

\end{theorem}

\section*{Competing interests}
\noindent The authors declare that they have no competing interests.

\vskip 0.3 cm

\section*{Acknowledgments}
\noindent The first author gratefully acknowledges that the ideas and concepts presented in this article were primarily motivated by the intellectual contributions and guidance of his late supervisor, Professor Dehici Abdelkader. This work continues his legacy in the field.


\end{document}